\documentclass[letterpaper, 10 pt, conference]{ieeeconf}
\IEEEoverridecommandlockouts

\usepackage{cite}
\usepackage{float}
\usepackage{amsmath,amssymb,amsfonts}
\usepackage{algorithmic}
\usepackage{mathrsfs}
\usepackage{graphicx}
\usepackage{textcomp}
\usepackage{xcolor}
\usepackage[amsmath,thmmarks]{ntheorem}
\usepackage{hyperref}
\usepackage{enumerate}
\hypersetup{pdfborder={0 0 0},}
\def\BibTeX{{\rm B\kern-.05em{\sc i\kern-.025em b}\kern-.08em
    T\kern-.1667em\lower.7ex\hbox{E}\kern-.125emX}}
\usepackage{algorithm}
\usepackage{algorithmic}
 \theoremheaderfont{\normalfont\bfseries} % 标题加粗
\theorembodyfont{\normalfont\itshape}    % 正文斜体
\theoremseparator{.}                     % 编号后添加句点

\newtheorem{remark}{Remark}
\newtheorem{theorem}{Theorem}
\newtheorem{lemma}{Lemma}

\newtheorem{corollary}{Corollary}

\newtheorem{assumption}{Assumption}

\newcommand{\AR}[1]{\textcolor{black}{{#1}}}

\begin{document}

\title{\LARGE \bf Decentralized Optimization via RC-ALADIN with Efficient Quantized Communication
}

\author{Xu Du,  Karl~H.~Johansson, and Apostolos I. Rikos$^*$
\thanks{$^*$Corresponding author.}
	\thanks{Xu Du and Apostolos I. Rikos are with the Artificial Intelligence Thrust of the Information Hub, The Hong Kong University of Science and Technology (Guangzhou), Guangzhou, China. 
    Apostolos I. Rikos is also affiliated with the Department of Computer Science and Engineering, The Hong Kong University of Science and Technology, Clear Water Bay, Hong Kong, China. E-mails: {\tt~michaelxudu@hkust-gz.edu.cn; apostolosr@hkust-gz.edu.cn}. 
            }
           % \thanks{Fumin Zhang is with the Department of Electronic and Computer Engineering, Hong Kong University of Science and Technology, Hong Kong SAR, China. E-mail: {\tt\small eefumin@ust.hk}.
            %} 
            \thanks{Karl H.~Johansson is with the Division of Decision and Control Systems, KTH Royal Institute of Technology, SE-100 44 Stockholm, Sweden. 
    He is also affiliated with Digital Futures, SE-100 44 Stockholm, Sweden. 
    E-mail:{\tt~kallej@kth.se}.
            }
            \thanks{The work of X.D. and A.I.R. is supported by the Guangzhou-HKUST(GZ) Joint Funding Scheme (Grant No. 2025A03J3960).}
           % \thanks{X.D. is also supported by the Guangzhou-HKUST(GZ) Joint Postdoctoral program (Grant No. P00560).}
}

\maketitle

\begin{abstract}
% \todo{Remove (Grant No. P00560) ??}
In this paper, we investigate the problem of decentralized consensus optimization over directed graphs with limited communication bandwidth. 
We introduce a novel decentralized optimization algorithm that combines the Reduced Consensus Augmented Lagrangian Alternating Direction Inexact Newton (RC-ALADIN) method with a finite time quantized coordination protocol, enabling quantized information exchange among nodes. 
% We introduce a novel decentralized optimization algorithm that integrates the Reduced Consensus Augmented Lagrangian Alternating Direction Inexact Newton (RC-ALADIN) method while enabling quantized information exchange among nodes. 
Assuming the nodes' local objective functions are $\mu$-strongly convex and simply smooth, we establish global convergence at a linear rate to a neighborhood of the optimal solution, with the neighborhood size determined by the quantization level. 
Additionally, we show that the same convergence result also holds for the case where the local objective functions are convex and $L$-smooth. 
Numerical experiments demonstrate that our proposed algorithm compares favorably against algorithms in the current literature while exhibiting communication efficient operation.

\end{abstract}

\section{Introduction}
% \todo{
% TO FIX: 
% % 1. Delete proofs (we will include them in the future version) (done)
% \\
% 2. 1st comment of 1st reviewer: 
% expand a bit explanations of assumptions. 
% \\
% Also, 4th comment of 1st reviewer: 
% explain importance of penalty parameter $\rho$, quantization level $\Delta$) (done)
% \\
% 3. 1st comment of 2nd reviewer: ``The authors should provide a more rigorous justification
% for adopting a two-layer algorithmic framework ...''(done)
% \\
% 4. 3rd comment from 3rd reviewer (done)
% \\
% 5. 1st comment of 4th reviewer (combined with comment 3.)
% \\
% 6. 3rd comment of 4th reviewer
% }
Decentralized optimization has garnered significant attention in recent years, driven by advances in areas such as control systems \cite{gosta2024decentralized}, machine learning \cite{decentralizedlearningsurvey}, and power grids \cite{liang2024decentralized}. This growing interest is largely due to the increasing need to address optimization problems that involve vast amounts of data and heterogeneous objective functions. 
% However, solving these problems in a centralized manner has been shown to be infeasible, as storing and processing such large data volumes on a single node is impractical \cite{boyd2011distributed}.

Decentralized optimization distributes data across multiple network nodes, typically using two main approaches: (i) primal decomposition and (ii) dual decomposition. Both methods locally optimize the objective functions.
Primal decomposition methods focus on sharing primal information among nodes, such as local optimal solutions or first- and second-order objective function data (e.g., \cite{nedic2009distributed, wang2018giant}). 
In contrast, dual decomposition methods update both the primal and the dual variables associated with coupling constraints (e.g., \cite{boyd2011distributed}, \cite{Houska2016}).
In general, dual decomposition methods achieve faster convergence and higher accuracy than primal decomposition approaches \cite[Section I]{ling2015dlm}. 
In this paper we focus on dual decomposition methods.

\noindent
\textbf{Existing Literature.} 
Dual decomposition-based optimization methods are built upon several core frameworks: Dual Decomposition \cite{dantzig1960decomposition}, the Alternating Direction Method of Multipliers (ADMM) \cite{boyd2011distributed}, and the Augmented Lagrangian Alternating Direction Inexact Newton (ALADIN) method \cite{Houska2016}. 
{\color{black}{While ADMM has undergone significant developments in both convergence theory and practical applications \cite{boyd2011distributed, he20121, Hong2016}, ALADIN improves convergence performance—achieving faster rates and ensuring convergence for both convex and nonconvex problems—by incorporating sequential quadratic programming techniques \cite{Houska2016, Du2023}. Subsequent developments have introduced specialized variants of ALADIN designed to tackle diverse computational complexity challenges \cite{Houska2021, Engelmann2019, Du2019, gosta2024decentralized, Du2023}.}}
However, two key limitations of these approaches are (i) their reliance on centralized coordination mechanisms (which limits scalability in distributed systems), and (ii) the requirement of nodes transmitting real-valued messages requiring a significant amount of bandwidth (which creates a scalability bottleneck). 
% This issue is compounded by two practical barriers: the need for full-duplex communication networks on undirected graphs, which is often impractical in real-world settings, and the reliance on real-valued information exchanges requiring infinite-bit precision, which creates a scalability bottleneck. These limitations become more pronounced as the network size grows, presenting a significant challenge in large-scale distributed optimization. 
To address the first challenge, the works \cite{falsone2023augmented, ling2015dlm, jiang2021asynchronous} investigated its decentralization for both resource allocation and consensus problems. 
Furthermore, recent works \cite{engelmann2020decomposition, gosta2024decentralized} have extended ALADIN to decentralized settings, however focusing only on resource allocation problems. 
Note that, despite some existing studies on variants of consensus ALADIN, such as \cite{Du2023, Du2023B}, a decentralized version of this framework has yet to be explored and addressed in the literature. 
Meanwhile, to address the second challenge of inefficient communication due to transmission of real-valued messages, the works in \cite{zhu2016quantized, elgabli2020q} investigated decentralized Consensus ADMM with quantized communication on undirected graphs. 
Although the aforementioned works have advanced ADMM development, they have only addressed some of the identified bottlenecks. 
The aforementioned literature has certain limitations, such as assuming a centralized coordinator, requiring nodes to exchange real-valued messages, or being restricted to undirected graphs. 
To the best of our knowledge, \cite{rikos2023asynchronous} is the only study to explore decentralized Consensus ADMM on directed graphs while enabling nodes to communicate in a resource efficient manner by exchanging quantized valued messages. 
Note that, although Consensus ALADIN exhibits superior convergence performance compared to Consensus ADMM \cite{Du2023}, its decentralized variant which however exhibits efficient communication with quantized value on directed graphs remains a problem that is unexplored in the literature. 

% \todo{example: 
% Motivated from the aforementioned challenges \\
% or gap in the literature \\
% or }

\noindent
\textbf{Main Contributions.} 
Motivated by the aforementioned challenges, we introduce a novel decentralized optimization algorithm leveraging Consensus ALADIN to tackle these issues. Note that to the authors knowledge, it is the first that Consensus ALADIN to address (i) fully decentralized algorithm operation, (ii) communication over directed graphs, and (iii) quantized communication among nodes. {\color{black}{Inspired by \cite{rikos2023asynchronous} and \cite{rikos2025distributed}, we note that adopting a two-layer algorithmic structure—where the optimization steps are decoupled from the decentralized averaging procedure—can achieve faster convergence compared to approaches that interleave optimization and communication in a single layer.}}
Our contributions are the following. 
\\ \noindent \textbf{A.} {\color{black}{We present a novel two-layer decentralized optimization algorithm, termed Quantized Decentralized Reduced Consensus ALADIN (QuDRC-ALADIN) (see Algorithm~\ref{alg:ALADIN}), designed to operate over directed communication graphs.
The inner layer (Algorithm~\ref{alg:QuAS}, inspired by \cite{rikos2023asynchronous}) enables efficient communication among nodes via the exchange of quantized messages.
The outer layer (i.e., the remaining steps of Algorithm~\ref{alg:ALADIN}) is responsible for updating the primal and dual variables of the original optimization problem.}}
%We present a novel decentralized optimization algorithm named Quantized Decentralized Reduced Consensus ALADIN (QuDRC-ALADIN) that operates over a directed communication graph. Our algorithm, enables nodes to communicate in an efficient manner by exchanging quantized valued messages (see Algorithm~\ref{alg:ALADIN}). 
\\ \noindent \textbf{B.} 
We prove that our algorithm converges to a neighborhood of the optimal solution (where the neighborhood depends on the quantization level) with global linear convergence rate for the case where the local cost function of each node is $\mu$-strongly convex, closed, proper, and simply smooth (see Theorem~\ref{the: convergence}). 
Additionally, we show that our algorithm exhibits global linear convergence rate also for the case where the local cost function of each node is convex, closed, proper, and $L$-smooth (see Corollary~\ref{coro1}).

\section{Notation and Preliminaries}\label{sec: Notation}

% \todo{change $a_i^k$ denotes the value of a vector $a$ of agent $i$ at iteration $k$ to $a_i^{[k]}$} 

\textbf{Notation.} 
We use the symbols $\mathbb{R}$, $\mathbb{Q}$, $\mathbb{Z}$, and $\mathbb{N}$ to represent the sets of real, rational, integer, and natural numbers, respectively. 
Matrices are indicated by capital letters (e.g., $A$), and vectors are represented by lowercase letters (e.g., $a$).
The transpose of matrix $A \in \mathbb R^{n \times n}$ and vector $a \in \mathbb R^{n}$ are represented as $A^\top$ and $a^\top$, respectively. 
For a real number $a \in \mathbb R $, $\lfloor a \rfloor$ and $\lceil a \rceil$ denote the greatest integer less than or equal to $a$ and the least integer greater than or equal to $a$, respectively. 
For the real vector $a \in \mathbb R^n$, $\lfloor a \rfloor \in \mathbb R^n$ and $\lceil a \rceil\in \mathbb R^n$ denote the element-wise operation. 
Furthermore, $\mathbf{1}$ represents the vector of all ones and $I$ denotes the identity matrix with appropriate dimensions. 
In addition, $\|a\|$ denotes the Euclidean norm of the vector $a$. 
The value of a variable $x$ of node $i$ at iteration $k$ is denoted as $x_i^{[k]}$.
The updated value is denoted as $(\cdot)^+$. 
%\todo{Furthermore, $\mathcal C$ stands for the closed, convex, proper set (CCP)},
Furthermore, $|\mathcal S|$ denotes the cardinality of a countable set $\mathcal S$ (e.g. $|\mathcal V|=N$ as we can see below). 
The notation $a|b$ denotes $b \in \mathbb R^n$ as the dual variable of constraint $a \in \mathbb R^n$. 

%\todo{C is defined below at Section III right?} 

%\todo{for sets and functions, proper has a different definition right? CCP can be defined only for functions the way we have it now. for sets proper means closed+convex+nonempty?} very important!

%\todo{I found that in books the term proper usually applies to functions. So we can call a function closed proper and convex. But it does not apply to sets. For sets, authors define closedness, convexity, and non-emptiness individually rather than a term like proper set. 
%}

%\todo{in my opinion the term proper for sets is not used so much, I think we could avoid.}

%\todo{lets discuss if we use the term CCP or simply closed, convex, proper (as it is done in books such as the book of Boyd: Distributed Optimization and Statistical
%Learning via the Alternating Direction
%Method of Multipliers}

\textbf{Graph Theory.} 
The communication network is captured by a directed graph (or digraph) $\mathcal G=(\mathcal V, \mathcal E)$. 
The set of agents (nodes) is denoted as $\mathcal V= \{1, \cdots, N\}$ (where $|\mathcal V|\geq 2$). 
The set of edges is denoted as $\mathcal E \subseteq \mathcal V\times \mathcal V \cup\{(i,i) \ | \ i \in \mathcal V \}$ (each node has a virtual self-edge).  
% With this definition,  each node $i$ also has an edge connected to itself. 
A directed edge from node $i$ to node $j$ is denoted by $e_{ji} \overset{\cdot}{=}(j,i)\in \mathcal E$. 
The subset of nodes that can directly transmit information to node $i$ is called the set of in-neighbors of $i$ and is denoted $\mathcal N_i^- = \{ j\in \mathcal V \ | \ (i,j)\in \mathcal E \}$. 
The subset of nodes that can directly receive information from node $i$ is called the set of out-neighbors of $i$ and is denoted $\mathcal N_i^+ = \{ l \in \mathcal V \ | \ (l,i) \in \mathcal E \}$. 
The cardinality of $\mathcal N_i^-$ represented as $\mathcal D_i^- = |\mathcal N_i^-|$, is called \emph{in-degree} of node $i$. 
The cardinality of $\mathcal N_i^+$ represented as $\mathcal D_i^+ = |\mathcal N_i^+|$, is called \emph{out-degree} of node $i$. 
The diameter $D$ of digraph $\mathcal G$ is the longest shortest path between $i , j \in \mathcal V$.  
A digraph $\mathcal G$ is strongly connected if there exists a directed path from every node $i$ to node $j$ that $i,j\in \mathcal V$. 

% Similarly, the \emph{out-neighbor} $\mathcal N_i^+=\{j\in \mathcal V| (i,j)\in \mathcal E \}$ denotes the subset of $\mathcal V$ where the information of node $i$ can be directly transmitted to. 
% Here, the cardinality of $\mathcal N_i^+$ is represented as $\mathcal D_i^+= |\mathcal N_i^+|$. 

\textbf{Quantization.} 
In digital communication networks, quantization serves to reduce bandwidth requirements and improve communication efficiency. 
By using a finite number of bits, quantization enables the application of error-correcting codes (e.g., Reed-Solomon, LDPC) to significantly enhance the signal's resilience to interference during transmission \cite{proakis2002communication}.
Three primary types of quantizers have been extensively studied: asymmetric, uniform, and logarithmic (see \cite{2019:Wei_Johansson}). 
In this paper we utilize asymmetric mid-rise quantizers with an infinite range (although our findings also hold for other types of quantizers as well). 
An asymmetric mid-rise quantizer is defined as 
\begin{equation}\label{eq: quantizer} %\small
    q^a_{\Delta}(b) = \left\lfloor \frac{b}{\Delta} \right\rfloor,
\end{equation}
where $b\in \mathbb R^n$ is the value to be quantized, and $\Delta \in \mathbb Q$ denotes the quantization level, and the superscript \emph{a} denotes the asymmetric type. 
% while $\left\lfloor \frac{b}{\Delta} \right\rfloor\in\mathbb Q$ are also rational number. 

% \ADD{APOSTOLOS: OK UNTIL HERE}

%\newpage
\section{Problem Formulation}\label{sec: problem}

%\REMOVE{In this section, we review the foundations of distributed and decentralized consensus optimization.} \todo{comment out?}
% and RC-ALADIN.

%\subsection{Distributed Optimization Problem Statement}\label{sec: DC}

Consider a communication network represented by a digraph $\mathcal{G} = (\mathcal{V}, \mathcal{E})$ comprising $N = |\mathcal{V}|$ nodes. 
We assume that the communication channels between nodes in our network $\mathcal{G}$ have limited bandwidth. 
Each node $i$ is associated with a scalar local cost function $f_i(x): \mathbb{R}^n \mapsto \mathbb{R}$, known exclusively to that node. 
Our objective in this paper is to develop a decentralized algorithm that enables nodes to collaboratively solve the following optimization problem 
\begin{equation}\label{eq: original} %\small
    \begin{split}
        \min_{z\in \mathbb R^n} & \quad \sum_{i=1}^{N} f_i(z),\quad i\in \{1,\cdots, N\},
    \end{split}
\end{equation}
where $z$ is the global optimization variable.
In order to solve problem \eqref{eq: original} %\REMOVE{by decomposing the local objectives,}  
we introduce a local variable $x_i$ for each node $i \in \mathcal{V}$ (following the approach outlined in \cite[Chapter 7.1]{boyd2011distributed}). 
Thus, \eqref{eq: original} is reformulated as 
\begin{equation}\label{eq: reformulate1}%\small
    \begin{split}
        \min_{x_i, \ i = 1, ..., N} & \quad \sum_{i=1}^{N} f_i(x_i)\\
        \text{s.t.}\;\;\quad& \quad x_i =z,\; \forall i\in\{1,\cdots, N\}, \quad x_i, z\in \mathbb R^n.
    \end{split}
\end{equation}
Problem \eqref{eq: reformulate1} is known as the  consensus optimization problem, as the constraint enforces equality among all local variables.
%\REMOVE{Our previous work RC-ALADIN \cite{Du2023}, which is based on Problem \eqref{eq: reformulate1}, is distributed with a centralized coordinator. 
%Specifically, while the updates of $x_i$s are performed in parallel, the update of 
%$z$ remains centralized.} %\todo{the last two sentences will go at the end of introducing ALADIN (Sec IV). We introduce ALADIN, we speak of its characteristics and then we speak of the shortcomings (e.g. the update of 
%$z$ remains centralized) and how we will solve them.}
% 
%\newpage
% 
%remove subsections
% 
%keep (3) and below (3) write (3) with constraeint of quantized comm. that new is prob (4). 
% 
%after new (4) we say ''The primary contribution of this paper...'' the black text
% 
%After Aladin, we say that z is the global (the way you have it). And because we do decentralized we introduce local copies etc. 
% 
% 
%\subsection{Modification of Distributed Optimization Problem}
% 
% In order to solve \eqref{eq: original} via the RC-ALADIN strategy \ADD{(see \cite{Du2023})} and guarantee efficient communication among nodes, we implement a key modifications. 
% First, we introduce the constraint $| x_i - x_j | \leq \epsilon$ for all pairs of nodes $i, j \in \mathcal{V}$, where $\epsilon \in \mathbb{R}$ is a predefined error tolerance. 
% Additionally, we impose the requirement that nodes exchange only quantized values. 
% With these adjustments incorporated, the problem \eqref{eq: original} is reformulated as 
In order to solve \eqref{eq: reformulate1} via the RC-ALADIN strategy (see \cite{Du2023}) while guaranteeing efficient communication among nodes, we have that \eqref{eq: reformulate1} is reformulated as 
\begin{equation}\label{eq: reformulate3}%\small
    \begin{split}
                \min_{x_i, \ i = 1, ..., N} & \quad \sum_{i=1}^{N} f_i(x_i)\\
        \text{s.t.}\;\;\quad& \;\quad x_i =z,\; \forall i\in\{1,\cdots, N\}, \quad x_i, z\in \mathbb R^n,\\
        &\;\quad \text{nodes communicate with quantized values.}
    \end{split}
\end{equation}
%The formulation of Problem \eqref{eq: reformulate3} also has the advantage of decoupling local problem-solving and information exchange into two distinct components: $x_i$ handles local optimization, while $z_i$ facilitates the interaction.

The primary contribution of this paper is the development of a decentralized algorithm that enables nodes to solve the problem in \eqref{eq: reformulate3}. 
In our paper, a decentralized approach refers to a distributed approach without the presence of a central coordinator. 
More specifically, in our proposed decentralized algorithm nodes coordinate solely through communication with their immediate neighbors (i.e., without the presence of a central coordinator). 
%In particular, the key contribution of this work is a fully decentralized algorithm for updating $z$ in \eqref{eq: reformulate3}. 
%\todo{the last line ''In particular, ...'' should go after presenting our algorithm. After we present our algorithm we will have a small paragraph comparing with previous works. In this paragraph we will say that ALADIN in [3] is distributed with a server node, and our current main contribution is to do it fully decentralized and in addition to guarantee efficient communication among nodes.}

\section{Preliminaries of RC-ALADIN}\label{sec: RC-ALADIN}

%\todo{I changed \cite{Du2023}.} 

% \todo{+ see my red comment from Section III-A. we need to mention at the end of IV, e introduce ALADIN, we speak of its
% characteristics and then we speak of the shortcomings (e.g.
% the update of z remains centralized) and how we will solve
% them.}

% \todo{1. before \eqref{eq: local update}, present Augmented Lagrangian}

% \todo{2. After Augmented Lagrangian, present \eqref{eq: RC-ALADIN} but before it speak of the local variables of each node (for example $x_i, g_i, \lambda_i$ and the global variable $z$ and what they are. For example $x_i$ is the optimal solution.}

The augmented Lagrangian of Problem \eqref{eq: reformulate3} is given by
\begin{equation}\label{eq: Lagrangian}%\small
   \mathscr L (x, z, \lambda) = \sum^N_{i=1} f_i(x_i) + \lambda_i^\top \left(x_i -z  \right) + \rho \left\| x_i -z \right\|^2,
 \end{equation}
where, $\rho>0$ is a penalty parameter, $\lambda=[\lambda_1^\top, \cdots,\lambda_N^\top]^\top$ denotes the dual variables and $x=[x_1^\top, \cdots,x_N^\top]^\top$ collects the local primal variables.
Focusing on \eqref{eq: Lagrangian}, RC-ALADIN was proposed in \cite{Du2023} to solve the consensus optimization problem in \eqref{eq: reformulate1}. 
Details of RC-ALADIN are the following,
\begin{subequations}\label{eq: RC-ALADIN}%\small
\begin{align}
        x_i^+ =& \mathop{\arg\min}_{x_i} f_i(x_i)+ \lambda_i^\top x_i + \frac{\rho}{2}\|x_i-z\|^2, \ \forall i \in \mathcal V , \label{eq: local update} \\
        g_i \ = &\rho\left(z-x_i^+\right)-\lambda_i, \ \forall i \in \mathcal V , \label{eq: gradient} \\
        z^+=&\frac{1}{N}\sum_{i=1}^N\left(x_i^+-\frac{g_i}{\rho}\right) , \label{eq: z update} \\
        \lambda_i^+=&\rho\left(x_i^+-z^+\right)-g_i, \ \forall i \in \mathcal V \label{eq: dual update}.  
\end{align}
\end{subequations}
In equation \eqref{eq: local update}, we minimize the augmented Lagrangian with respect to each $x_i$. In equation \eqref{eq: gradient}, we evaluate the (sub)gradient of each $f_i$ at $x_i^+$. 
Note that \eqref{eq: gradient} always exists as long as \eqref{eq: local update} is solvable. 
Equations \eqref{eq: z update} and \eqref{eq: dual update} are the closed form expression of the following centralized reduced consensus QP problem 
\begin{equation}\label{eq: reduced QP} %\small
    \begin{split} 
        \min_{\Delta x_i \in \mathbb R^n, \ i = 1, ..., N} & \quad \sum_{i=1}^{N} \frac{\rho}{2}\Delta x_i^\top\Delta x_i+g_i^\top \Delta x_i\\ 
        \text{s.t.}\hspace{2mm}\qquad& \;\quad x_i^+ +\Delta x_i =z |\lambda_i,\;\forall i\in\{1,\cdots, N\}.
    \end{split} 
\end{equation} 
%\todo{is min correct? it has $z \in \mathcal{R}$ below it.}
The reduced consensus QP in \eqref{eq: reduced QP} is used from nodes to coordinate without depending on second-order information of their local cost functions $f_i$, $\forall i \in \mathcal V$. 
Our previous work, RC-ALADIN in \eqref{eq: RC-ALADIN}, provides global convergence guarantees for \eqref{eq: reformulate1} when the local functions of each node are convex. Additionally, it offers globally linear convergence guarantees for \eqref{eq: reformulate1} when the local functions of each node are $\mu$-strongly convex (see \cite[Theorem~$2$, Theorem~$7$]{Du2023_Arxiv}).

Note that, RC-ALADIN in \eqref{eq: RC-ALADIN} is implemented by nodes in a distributed manner, but also considers the presence of a centralized coordinator that can exchange messages with every node in the network. 
More specifically, focusing on \eqref{eq: local update} -- \eqref{eq: dual update}, the updates \eqref{eq: local update}, \eqref{eq: gradient}, \eqref{eq: dual update} are performed locally from each node $i$, and the update \eqref{eq: z update} is performed in a centralized fashion from the centralized coordinator. 
%\todo{Specifically, while the updates of $x_i$s are performed in parallel (as shown in \eqref{eq: local update}), the update of $z$ remains centralized (as described in \eqref{eq: z update}).} 
%\todo{A fully decentralized variant of RC-ALADIN will be introduced in the next section, where the centralized update in \eqref{eq: z update} is replaced by a decentralized operation through the introduction of local copies $z_i$s of $z$.} 
Motivated by this limitation, in the next section we will present a fully decentralized algorithm. 
Our proposed algorithm enables nodes to (i) collaboratively  solve problem \eqref{eq: reformulate3} by communicating exclusively with their immediate neighbors (eliminating the need for a central coordinator), and (ii) exhibit efficient communication within the network. 

% \todo{
% \\ this part above: ''the centralized update in \eqref{eq: z update} is replaced by a decentralized operation through the introduction of local copies $z_i$s of $z$.''  will go to the comparison paragraph I spoke in Section III.  
% I wrote: \\ 
% In this paragraph we will say that ALADIN in [3] is distributed with a server node, and our current main contribution is to do it fully decentralized and in addition to guarantee efficient communication among nodes}

\section{Decentralized RC-ALADIN with Efficient Communication}\label{sec: algorithm}

%\todo{BE CONSISTENT WITH NOTATION. WHAT IS $\hat z_i$? IS IT global variable estimation? OR local copies? WE SHOULD USE ONE NAME. I think we should use ``global variable estimation''. please fix this in all the paper.} 

%\todo{lets avoid the work variation of [3]. It is a new work inspired from [3]}. 

In this section, we introduce a novel decentralized algorithm designed to address the problem in \eqref{eq: reformulate3}. 
%This algorithm draws inspiration from RC-ALADIN \eqref{eq: RC-ALADIN}. \todo{comment out last sentence. After introducing our algorithm, we will add a paragraph speaking of novelties and differences with previous works. In this paragraph we will say that it is inspired from RC-ALADIN \eqref{eq: RC-ALADIN}.} 
Before introducing the proposed decentralized algorithm, we first make the following assumptions that are important for our subsequent development. 
\begin{assumption}\label{ass:1}
  %  \todo{Assume that $\mathcal G$ is strongly connected.} 
    The communication network is modeled as a \textit{strongly connected} digraph $\mathcal{G} = (\mathcal{V}, \mathcal{E})$. 
    Also, every node $i$ knows the diameter of the network $D$, and a common quantization level $\Delta$. %\todo{please fix explanation below since we updated this assumption}
    %and a finite number of upper bounds $\mathcal B\in \mathbb N$ for the nodes gets consensus (average). 
\end{assumption}
% \begin{assumption}\label{ass:3} \todo{keep this for later. after proof we can say that if Assump3 does not hold and we only have Assump2 then we can have a rate of convergence that is ... this}
%     The local cost function $f_i$ of each node $i \in \mathcal{V}$ is smooth, closed, proper and convex. 
%     There is a unique global minimizer of \eqref{eq: reformulate3}. \todo{last sentence is the explanation for strong convexity and smoothness right? it needs to be mentioned below. For example: Assumption~\ref{ass:3} ensures that there is a unique global minimizer of \eqref{eq: reformulate3}.} 
%     \todo{Also, ''smooth'' is what we present in Assumption~\ref{ass:2}? why mention it again?} 
% \end{assumption} 

\begin{assumption}\label{ass:2} 
The local cost function $f_i$ of each node $i \in \mathcal{V}$ is closed, proper, simply smooth and $\mu$-strongly convex. 
Specifically, for each local cost function $f_i$, for every $x_\alpha,x_\beta \in \mathbb R^n$, there exist a \emph{strong-convexity} constant $\mu_i>0$ such that 
\begin{equation}\label{eq: mu}
\begin{split}
    f_i(x_\alpha)+&\nabla f_i(x_\alpha)^\top (x_\beta-x_\alpha) +\frac{\mu_i}{2}\|x_\beta-x_\alpha\|^2\leq f_i(x_\beta). 
\end{split}
\end{equation}
Here, simply smooth refers to a smooth function $f_i$  where $L_i$ cannot be explicitly estimated as in inequality \eqref{eq: lip} in Assumption~\ref{ass:3} below. 
In contrast, $L$-smooth indicates that $f_i$ is smooth with an explicitly estimable $L_i$ as shown in \eqref{eq: lip} below.
\end{assumption}

\begin{assumption}\label{ass:3} 
The local cost function $f_i$ of each node $i \in \mathcal{V}$ is closed, proper, $L$-smooth and  convex. 
Specifically, for each local cost function $f_i$, for every $x_\alpha,x_\beta \in \mathbb R^n$, there exist a \emph{Lipschitz-continuity} constant $L_i>0$ (see \cite[equation (51)]{ling2015dlm}) such that 
\begin{equation}\label{eq: lip}%\small
\begin{split}  
     &  \left\|\nabla f_i(x_\alpha)-\nabla f_i(x_\beta) \right\|^2 \\
     &\leq L_i \left(\nabla f_i(x_\alpha)-\nabla f_i(x_\beta)\right)^\top\left(x_\alpha-x_\beta \right).
\end{split}
\end{equation}
\end{assumption}

% \todo{TWO ASSUMPTIONS: ASsum2: SIMPLY SMOOTH AND M-STRONGLY CONVEX (and also in assump2 we mention that simply smooth means that it is smooth but we dont know what is the parameter L etc). \\ done
% Assump3: L-SMOOTH AND CONVEX. \\ 
% mention that simply smooth means it is smooth but we dont know the L parameter
% L-smooth means that it is smooth and we know the L parameter.\\ done
% so in paper we use ``simply smooth'' and ``L-smooth'' !! THIS IS IMPORTANT \\done
% IN theorem 1 we mention that we have Assumption 2 and we prove convergence. 
% In Corollary 1 we mention we have Assumption 3 and we prove convergence. 
% }

%\todo{why there is a . after the Assumption 1 here? We need to find a way to delete the dots. same for assumption 2}
In Assumption~\ref{ass:1}, strong connectivity ensures that information can propagate between all nodes in the network and guarantees the convergence of Algorithm~\ref{alg:QuAS} (since strongly connected digraph implies there is a path from every node to every other node in the network).  
Knowing the diameter of the network is useful for each node to determine whether Algorithm \ref{alg:QuAS} has converged, allowing it to proceed to step~$3$ of Algorithm~\ref{alg:ALADIN}. 
\AR{
Note here that the network diameter can can be computed by employing distributed algorithms \cite{OLIVA201620}. 
}
{\color{black}{
For open undirected networks, the network diameter 
$D$ need not be known a priori \cite{deplano2025optimization}. 
}}
The quantization level is important for quantizing the messages as described in \eqref{alg:QuAS} (thus ensuring efficient communication during the execution of Algorithm \ref{alg:QuAS}). 
\AR{Nodes can compute a common quantization level $\Delta$ in finite time via a max-consensus operation.}
Assumption~\ref{ass:2} ensures strong convexity and simple smoothness, allowing us to establish a global linear convergence rate for Algorithm~\ref{alg:ALADIN} while ensuring that the global cost function in \eqref{eq: reformulate3} has a unique minimum \cite[Theorem~13.27]{beck2017first}, \cite{ling2015dlm}. 
Assumption~\ref{ass:3} guarantees that the Lipschitz continuity of gradients in \eqref{eq: lip} ensures the existence of a global optimal solution $x^*$ for \eqref{eq: reformulate3} and enables nodes to compute it, which is a standard requirement in first-order distributed optimization frameworks (see, e.g., \cite{beck2017first}). {\color{black}{Many practical optimization problems satisfy Assumptions~\ref{ass:2} and \ref{ass:3} \cite{Rawlings2017,boyd2004convex}. 
However, when these two assumptions are violated in specific applications, the proposed Algorithm~\ref{alg:ALADIN} admits only global convergence, without a guaranteed convergence rate \cite{Du2023_Arxiv,Du2023}.}}

\subsection{Algorithm Development}\label{sec: algorithm structure}

In this section, we present our proposed decentralized algorithm, 
detailed below as Algorithm~\ref{alg:ALADIN}.

\begin{algorithm}[h] 
	%\small
	\caption{QuDRC-ALADIN: Quantized Decentralized Reduced Consensus ALADIN}
	\textbf{Input.} Strongly connected digraph $\mathcal{G} = (\mathcal{V}, \mathcal{E})$, parameter $\rho$, network diameter $D$, quantization level $\Delta$, for each node $i \in \mathcal{V}$. 
    Each node $i \in \mathcal V$ has a local cost function $f_i$. 
    Assumptions~\ref{ass:1} and \ref{ass:2} hold. 
    \\
    \textbf{Initialization.} Randomly chosen dual variable $\hat \lambda_i \in \mathbb{R}^n$, and global variable estimation $\hat z_i \in \mathbb{R}^n$, for each node $i \in \mathcal{V}$. 
    \\
    \textbf{Iteration.} 
	Each node $i \in \mathcal{V}$ repeats:
	\begin{enumerate}
	\item  Optimize $x_i$ as 
    \begin{equation}\label{eq: new local primal}
        x_i^+ = \mathop{\arg\min}_{x_i} f_i(x_i)+ \hat\lambda_i^\top x_i + \frac{\rho}{2}\|x_i-\hat z_i\|^2.
    \end{equation}
    \item Evaluate the gradient $g_i$ of $f_i$ as 
    \begin{equation}\label{eq: new gradient}
         g_i = \rho\left(\hat z_i-x_i^+\right)-\hat\lambda_i.
    \end{equation}
    \item  Calculate the global variable estimation as  
    \begin{equation}\label{eq: local copy update}
        \hat z_i^+= \text{Algorithm~\ref{alg:QuAS}} \left(x_i^+-\frac{g_i}{\rho}, D, \Delta \right).
    \end{equation}
    \item Update the dual variable $\hat\lambda_i$ as 
    \begin{equation}\label{eq: new dual}
        \hat\lambda_i^+=\rho\left(x_i^+-\hat z_i^+\right)-g_i.
    \end{equation}
	\end{enumerate}
    \textbf{Output.} Each node $i$ calculates $x_i^*$ that solves problem \eqref{eq: reformulate3}. 
	\label{alg:ALADIN}
\end{algorithm}
The intuition of Algorithm~\ref{alg:ALADIN} is organized into two main phases: local optimization and coordination among nodes. 
In the first step, each node $i \in \mathcal{V}$ performs a local optimization to determine the optimal value of its variable $x_i$ by solving its corresponding augmented objective function, (see \eqref{eq: new local primal}). 
Following this, in the second step each node $i$ evaluates the (sub)gradient of its local function $f_i$ at the locally optimized solution $x_i^+$. 
This (sub)gradient evaluation serves as preparation for the aggregation process in the subsequent step (see \eqref{eq: new gradient}). 
In the third step, all nodes collaborate to update their estimates of the global variable $\hat z_i^+$ through the quantized, decentralized operation of Algorithm~\ref{alg:QuAS} (see \eqref{eq: local copy update}). 
Finally, in the fourth step each node updates its dual variable $\hat\lambda_i^+$ which encodes sensitivity information related to the constraints of problem~\eqref{eq: reformulate3}. 
The updated dual variables are then utilized in the next iteration's local optimization phase (see \eqref{eq: new dual}). 
Overall, Algorithm~\ref{alg:ALADIN} alternates between performing local optimization (step~$1$) and solving problem~\eqref{eq: reduced QP}, iterating until convergence is achieved and the optimal solution is obtained.
 
\begin{algorithm}[h]
	% \small
	\caption{FQAC: Finite-time Quantized Average Consensus}
	\textbf{Input.} $y_i=x_i^+-\frac{g_i}{\rho}, D, \Delta$. \\
    \textbf{Initialization.} Each node $i \in \mathcal{V}:$
    \begin{enumerate}
	\item Assigns probability 
    \begin{equation}%\small 
    p_{li}=\left\{
   \begin{split}
       &\frac{1}{1+\mathcal D_i^+},\quad& \text{if}\; l\in \mathcal N_i^+\cup \{i\},\\
       &0,\quad & \text{if}\; l\notin \mathcal N_i^+\cup \{i\},
   \end{split}
      \right.   
    \end{equation} 
    to each out-neighbor of node $i$.
    \item Sets $\xi_i = 2$, $\chi_i= 2 q^a_{\Delta}(y_i)$ (see \eqref{eq: quantizer}).
	\end{enumerate}
	\textbf{Iteration.} For time steps $t=1,2,\cdots$ each node $i \in \mathcal V$ does: 
	\begin{enumerate}
	\item \textbf{If} $t\;\text{mod}(D) = 1$, sets $M_i = \left\lceil \frac{\chi_i}{\xi_i} \right\rceil$ and $m_i=\left\lfloor \frac{\chi_i}{\xi_i} \right\rfloor$.\vspace{2mm}

\item %\todo{IMPORTANT: When we have in-neighbor use $l$, when we have out-neighbor use $j$} \\
Broadcasts $M_i, m_i$ to each out-neighbor $l \in \mathcal N_i^+$ and receives $M_j, m_j$ from each in-neighbor $j \in \mathcal N_i^-$. 
Then, sets $M_i = \text{max}_{j \in \mathcal N_i^- \cup \{i\}}\; M_j$, $m_i = \text{min}_{j \in \mathcal N_i^-\cup \{i\}}\; m_j$. \vspace{2mm}
\item Sets $\tau_i = \xi_i$. \vspace{2mm}
\item \textbf{While} $\tau_i>1$ \textbf{do}
\begin{enumerate}
    \item $c_i=\left\lfloor \frac{\chi_i}{\xi_i} \right\rfloor$. \vspace{2mm}
    \item Sets $\chi_i=\chi_i - c_i$, $\xi_i=\xi_i-1$, $\tau_i=\tau-1$. 
    \item Transmits $c_i$ to randomly chosen out-neighbor $l \in \mathcal N_i^+\cup \{i\}$ with probability $p_{li}$. 
    \item Receives $c_i$ from $j\in \mathcal N_i^-$ and updates
    \begin{subequations}\label{eq: local information avg}%\small
    \begin{align}
    \chi_i^{[t+1]} &= \chi_i^{[t]} + \sum_{j\in \mathcal N_i^- }w_{ij}^{[t]} c_j^{[t]}, \\
    \xi_i^{[t+1]} &= \xi_i^{[t]} + \sum_{j\in \mathcal N_i^-} w_{ij}^{[t]}.
        \end{align}
    \end{subequations}    
    Here $w_{ij}^{[t]} = 1$ if node $i$ receives $c_j^{[t]}$ from node $j$ at step $t$. 
    Otherwise $w_{ij}^{[t]}=0$ and node $i$ does not receive information from node $j$. 
\end{enumerate}
\item \textbf{if} $t\; \text{mod}\; (D)=0$ and $\left\|M_i-m_i\right\|_{\infty}\leq 1$, set $\hat z_i^+ = m_i \Delta$, and stop the operation of the algorithm. 
	\end{enumerate}
    \textbf{Output.} $\hat z_i^+$.
	\label{alg:QuAS}
\end{algorithm}

Algorithm~\ref{alg:QuAS} follows a structure similar to \cite[Algorithm~$1$]{rikos2022non}, and consists of three main operations: quantization, averaging, and a stopping criterion.  
During initialization each node $i$ quantizes its local information $y_i = x_i^+ - \frac{g_i}{\rho}$ into a quantized value $\chi_i$. 
Then, it splits $\chi_i$ into $\xi_i$ pieces (the value of some pieces might be greater than others by one). 
It retains the piece with the smallest value to itself and transmits the rest $\xi_i - 1$ pieces to randomly chosen out-neighbors $l \in \mathcal N_i^+$ or to itself. 
Then, it receives the pieces $c_j$ transmitted from each in-neighbor $j \in \mathcal N_i^-$ and updates $\chi_i$ and $\xi_i$ as in \eqref{eq: local information avg}. 
The algorithm also performs max- and min-consensus operations every $D$ time steps. 
If the results of the max-consensus $M_i$ and min-consensus $m_i$ have a difference less or equal to one, then each node $i$ scales its solution according to the quantization level to compute $\hat z_i^+$.
At this point, Algorithm~\ref{alg:QuAS} terminates, and each node $i$ transitions to step~$4$ of Algorithm~\ref{alg:ALADIN}. 
{\color{black}{Note that Algorithm~\ref{alg:QuAS} is guaranteed to converge in finite time, ensuring $\left\|M_i-m_i\right\|_{\infty}\leq 1$ (see Step 5). 
The convergence time depends on the network diameter $D$, as established in \cite[Theorem~1]{rikos2024distributed}.}}

% \todo{DELETE BELOW. WAS REPLACED WITH ABOVE}
% Algorithm \ref{alg:QuAS} has a similar structure as \cite[Algorithm 1]{rikos2022non}. 
% It contains three main operations (i) quantization, (ii) averaging, and (iii) a stopping criterion. 
% During each iteration $t$ of Algorithm \ref{alg:QuAS}, each node $i$ quantizes its local information $y_i =x_i^+-\frac{g_i}{\rho}$ as $\chi_i$. 
% If the token $\tau_i$ of node $i$ is positive, then:
% \begin{itemize}
%     \item splits $\chi_i$ into $\xi$ pieces (the value of some pieces might be greater than others by one);
% \item transmits each piece to a randomly chosen out-neighbor $l \in \mathcal N_i^+$ including itself;
% \item receives pieces $c_j$ from its in-neighbor $j\in \mathcal N_i^-$, updates $\chi_i$ and $\xi_i$ locally as equation \eqref{eq: local information avg}.
%     \end{itemize}
%  %   \footnotetext[1]{Here $w_{ij}^{[t]} = 1$ if node $i$ receives $c_j^{[t]}$ from node $j$ at step $t$. Otherwise set $w_{ij}^{[t]}=0$ and node $i$ do not receive information from node $j$. 
% %}
% Finally, if the results of the max-consensus $M_i$ and
% min-consensus $m_i$ have a difference less or equal to one, each node $i$ scales the solution according to the quantization level as $\hat z_i^+$ and continue at step~$4$ of Algorithm~\ref{alg:ALADIN}. 
% Note that Algorithm~\ref{alg:QuAS} is guaranteed to converge in a finite time \cite{rikos2022non}.

\vspace{.2cm}
\noindent
\textbf{Comparison with Previous Works.} 
Note that the key contribution of this work is the development of a communication-efficient fully decentralized algorithm for solving problem~\eqref{eq: reformulate3}.
More specifically, while our previous work RC-ALADIN in \eqref{eq: RC-ALADIN} (also see \cite{Du2023}) is distributed and considers the existence of a server node able to communicate with every node in the network, our main contribution in this paper lies in achieving full decentralization (i.e., nodes coordinate without the presence of a server node). 
This is achieved by replacing the centralized update in \eqref{eq: z update} by a decentralized operation through the introduction of local copies $\hat z_i$ of $z$. 
Moreover, by implementing Algorithm~\ref{alg:QuAS}, Algorithm~\ref{alg:ALADIN} enhances its communication-efficiency while ensuring convergence precision (see for example \cite{rikos2023asynchronous, rikos2023distributed, rikos2024distributed}). 
This characteristic is not present in our previous work RC-ALADIN in \cite{Du2023} where nodes are operating in a resource inefficient manner by exchanging real-valued messages that require a significant amount of bandwidth.  
Additionally, in our previous work \cite{rikos2023asynchronous}, the convexity (without simple smoothness) of $f_i$ for all $i$ is sufficient for establishing convergence, demonstrating a global \textit{sub-linear} convergence rate. 
Furthermore, the algorithms in \cite{rikos2023distributed, rikos2024distributed} require $L$-smoothness \textit{and} $\mu$-strong convexity to establish the global \textit{linear} convergence rate. 
In contrast, in this paper, \textit{either} simple smoothness with $\mu$-strong convexity \textit{or} $L$-smoothness with simple convexity of $f_i$ for all $i$ is required to establish the global linear convergence rate of Algorithm~\ref{alg:ALADIN}. 
Finally, note that our previous works \cite{Du2023, Du2023_Arxiv} require the same assumptions as this paper (i.e., \textit{either} simple smoothness with $\mu$-strong convexity \textit{or} $L$-smoothness with simple convexity of $f_i$ for all $i$). 
However, as we mentioned above they do not exhibit communication efficient operation among nodes. 
Details of this will be provided in the next subsection.

\subsection{Convergence Analysis}\label{sec: convergence}

In this section, we provide the convergence analysis of Algorithm~\ref{alg:ALADIN}. 
First, we introduce the two lemmas that are important for our analysis. 
Then, we prove our main result via a theorem. 
%Note that our theorem is inspired by \cite[Theorem~$7$]{Du2023_Arxiv}. 
%\todo{not sure last sentence is needed. will check again}

%\todo{FOR THE LEMMA BELOW PLEASE CHECK: \cite[Lemma~$1$]{rikos2023distributed2}.} 

\begin{lemma} 
The update of the global variable estimation $\hat z_i$ of each node $i \in \mathcal{V}$ is given by Algorithm~\ref{alg:ALADIN} in \eqref{eq: local copy update}. 
% In this update an error $\epsilon$ related to the quantization level $\Delta$ is introduced. 
% Note that for the error $\epsilon$ it holds that $\epsilon \geq 2 \Delta$ \todo{not $\epsilon \leq 2 \Delta$ ???}. 
According to the constraints of problem~\eqref{eq: reformulate3}, the following equation is satisfied 
% \begin{equation}\label{eq: range} 
%         \hat z_i^+\in \left[\frac{1}{N} \sum_{i=1}^N y_i -2\Delta,\frac{1}{N} \sum_{i=1}^N y_i +2\Delta\right], 
%     \end{equation} 
%     where $y_i = x_i^+-\frac{g_i}{\rho}$. 
    %This means 
    \begin{equation}\label{eq: error}\small
    \left\{
    \begin{split} 
        &\hat z_i^+ = \frac{1}{N} \sum_{i=1}^N \Delta \left \lfloor \frac{y_i}{\Delta} \right \rfloor + \kappa_i, \;\left\| \kappa_i \right\|_{\infty}\leq \Delta, \\
        &\left\|z^+-\hat z_i^+\right\|_\infty\leq 2\Delta,\\[2mm]
                %\hat \lambda_i^+ - \lambda_i^+ &= \rho(z^+-\hat z_i^+),  \; 
                &\left\| \hat \lambda_i^+ - \lambda_i^+
 \right\|_{\infty}= \left\|\rho(z^+-\hat z_i^+)\right\|_{\infty}\leq 2\rho \Delta,
    \end{split}\right. 
\end{equation}
where $y_i = x_i^+-\frac{g_i}{\rho}$.

% \todo{With the local copy update $\hat z_i$ given by the quantized averaging algorithm in \eqref{eq: local copy update}, an error $\epsilon$ related to the quantization level $\Delta$ is introduced. 
% It is important to note that $\epsilon \geq 2 \Delta$. 
% According to the constraints of Problem \eqref{eq: reformulate3}, the following equation will be satisfied,
% \begin{equation}\label{eq: range}%\small
%         \hat z_i^+\in \left[\frac{1}{N} \sum_{i=1}^N y_i -\frac{\epsilon}{2},\frac{1}{N} \sum_{i=1}^N y_i +\frac{\epsilon}{2}\right].
%     \end{equation}
%     This means \begin{equation}\label{eq: error}\left\{
%     \begin{split}
%         \hat z_i^+ &= \frac{1}{N} \sum_{i=1}^N \Delta \left \lfloor \frac{y_i}{\Delta} \right \rfloor + \kappa_i, \;\left\| \kappa_i \right\|_{\infty}\leq \Delta;\\
%                 \hat \lambda_i^+ - \lambda_i^+ &= \rho(z^+-\hat z_i^+),  \; \left\| \hat \lambda_i^+ - \lambda_i^+
%  \right\|_{\infty}\leq \rho \Delta.
%     \end{split}\right.
% \end{equation}} 
\end{lemma}
\textit{Proof.} See \cite[Lemma~$1$]{rikos2023distributed2}.
% \cite[Remark~$1$, equation~$(24)$]{rikos2023asynchronous}, \cite[equation~$(11)$]{rikos2023distributed2} and  
%  \cite[Lemma~$1$]{rikos2023distributed2}.     \todo{lets check again and focus on \cite[Lemma~$1$]{rikos2023distributed2}}  
 \hfill$\blacksquare$

\begin{lemma}
% For the distributed consensus \ADD{optimization} problem in~\eqref{eq: reformulate3}, with Algorithm \ref{alg:ALADIN}, the local primal update $x_i^+$ has a relationship with the local dual $\hat \lambda_i$, $\hat \lambda_i^+$ and global primal variables approximation $\hat z_i$, $\hat z_i^+$ in the following way 
For the distributed consensus optimization problem presented in \eqref{eq: reformulate3}, Algorithm \ref{alg:ALADIN} establishes a relationship between the local primal update $x_i^+$, the local dual variables $\hat \lambda_i$ and $\hat \lambda_i^+$, and the global primal variable approximations $\hat z_i$ and $\hat z_i^+$. 
This relationship is 
    \begin{equation}\label{eq: lemma}%\small
        x_i^+ = \frac{\hat \lambda_i^+-\hat\lambda_i}{2\rho}+\frac{\hat z^+_i+\hat z_i}{2}.
    \end{equation}
\end{lemma}
\textit{Proof.} From \eqref{eq: new dual} we have 
\begin{equation}\label{eq: lemma2}%\small
    \begin{split}
       \hat \lambda_i^+&=\rho\left(x_i^+-\hat z_i^+\right)-g_i\\
        &\overset{\eqref{eq: gradient}}{=} \rho\left(x_i^+-\hat z_i^+\right) -  \left(\rho\left(\hat z_i-x_i^+\right)-\lambda_i \right)\\
        &= 2\rho x_i^+ - \rho(\hat z_i^++\hat z_i) +\hat\lambda_i.
    \end{split}
\end{equation}
equation \eqref{eq: lemma} is then derived from \eqref{eq: lemma2}. \hfill$\blacksquare$

Moreover, from the KKT (Karush-Kuhn-Tucker) system of problem \eqref{eq: reduced QP} and \eqref{eq: error}, the following formulas can be obtained,
\begin{equation}\label{eq: lam}\small
    \sum_{i=1}^N \lambda_i^{+} = 0, \ \text{and} \ \left\|\sum_{i=1}^N \hat\lambda_i^{+}\right\|_\infty \leq 2\rho N \Delta.%\sum_{i=1}^{N} \kappa_i.
\end{equation}
In \eqref{eq: lam} the first equality arises from the KKT stationarity condition for the reduced QP problem in \eqref{eq: reduced QP}. 
The second inequality arises from substituting the third inequality of \eqref{eq: error} into \eqref{eq: new dual} and then summing over all nodes. 
% \todo{IS THIS CORRECT? also is it $2\rho N \Delta$ or $\rho N \Delta$ ? because $\kappa_i \leq \Delta$... I think it is $\rho N \Delta$ in the second equality above... also is it $= \rho N \Delta$? or $\leq \rho N \Delta$? because $\kappa_i \leq \Delta$} 
For establishing the global convergence of Algorithm~\ref{alg:ALADIN}, we introduce the following Lyapunov function, %\cite[Section 5.1]{Du2023_Arxiv} \todo{no ref needed here I think},
\begin{equation}\label{eq: lya}\small
    \mathcal L(z,\lambda) = \frac{1}{\rho} \sum_{i=1}^N \left\|\lambda_i - \lambda_i^* \right\|^2 + \rho N\left\|z - z^* \right\|^2,
\end{equation}
where $\lambda=[\lambda_1^\top, \lambda_2^\top, \cdots, \lambda_N^\top]^\top$, and $(z^*, \lambda^*)$ denotes the optimal solution pair of problem \eqref{eq: reformulate3}.
Following the analysis in \cite{jiang2021asynchronous}, we define the finite positive scalars $0<M_z< \infty$  for simplifying our later analysis, such that  $\left\|z -z^* \right\|\leq M_{z}$
%and $\left\|\hat \lambda_i -\lambda_i^* \right\|\leq M_{\lambda}$, 
for every node $i \in \mathcal{V}$. 
\AR{Note that the proof of the theorem below will be available in an extended version of our paper.}

% \todo{above, we have $\left\|\hat z_i - z^* \right\|\leq M_{z}$, $\left\|z -z^* \right\|\leq M_{z}$ ... they are bounded from the same $M_{z}$?}

\begin{theorem}\label{the: convergence}
Let us consider a digraph $\mathcal G= \left( \mathcal V, \mathcal E\right)$. 
Each node $i \in \mathcal V$ has a local cost function $f_i$, and Assumptions~\ref{ass:1} and \ref{ass:2} hold. 
Each node $i \in \mathcal V$ in the network executes Algorithm~\ref{alg:ALADIN} for solving the consensus optimization problem in~\eqref{eq: reformulate3} in a decentralized fashion. 
Given parameter $\rho>0$, % and an error tolerance $\epsilon \in \mathbb Q$, %\todo{(here we say parameter $\rho>0$,  error tolerance $\epsilon$. .. we need to mention the parameters we are using. if we dont use D then we dont need to mention it.)... in \eqref{eq: condition} we use only $\delta$, and after it we use $\rho$. why we mentioned D?}
during the operation of Algorithm~\ref{alg:ALADIN} there always exists a $\delta>0$ such that 
\begin{equation}\label{eq: condition}%\small
    \delta \mathcal{L}\left(\hat z^+, \hat\lambda^+\right)\leq 4 \sum_{i=1}^N \mu_i\left\|x_i^+ -z^* \right\|^2 ,  
\end{equation} 
where $\hat z = \hat z_i, \forall i \in \mathcal{V}$. 
From \eqref{eq: condition}, we have that during the operation of Algorithm~\ref{alg:ALADIN} the following inequality is satisfied 
\begin{equation}\label{eq: reault1}%\small
    \mathcal L \left(\hat z^+, \hat\lambda^+\right) \leq \frac{1}{1+ \delta}\mathcal L \left(\hat z, \hat\lambda\right) + \frac{4}{1+ \delta} \mathcal O(N\Delta),
\end{equation} 
where $\mathcal O(N \Delta) = 6\rho  M_z  N\Delta$,
% \begin{equation}
%     \mathcal O(N \Delta) = 6\rho  M_z  N\Delta, 
% \end{equation}
 $\Delta$ is the utilized quantization level, and $\left\|z -z^* \right\|\leq M_{z} <\infty$. 
% For the distributed consensus optimization problem in \eqref{eq: reformulate3} defined on a strongly connected digraph $\mathcal G= \left( \mathcal V, \mathcal E\right)$, let Assumption~\ref{ass:1} and Assumption~\ref{ass:2}  be satisfied for the local cost function $f_i: \mathbb R^n\rightarrow \mathbb R$ of each node $i\in \mathcal V$. 
% Let $f_i$s be $\mu_i$ strongly convex \eqref{eq: mu}. \todo{we said that Assumption-2 holds already.}
% Given parameter $\rho>0$, network diameter $D$,  an error tolerance $\epsilon \in \mathbb Q$.
%     %Define $\bar z = \frac{1}{N} \sum_{i=1}^N \hat z_i$, $x=[x_1^\top, x_2^\top, \cdots, x_N^\top]^\top$, $\lambda=[\lambda_1^\top, \lambda_2^\top, \cdots, \lambda_N^\top]^\top$. 
% During the operation of Algorithm~\ref{alg:ALADIN}, there always exists a $\delta>0$ such that
% \begin{equation}\label{eq: condition}
%     \delta \mathcal{L}(\hat z^+, \hat \lambda^+)\leq 4 \sum_{i=1}^N \mu_i\left\|x_i^+ -z^* \right\|^2,
% \end{equation}
% where $\hat z = \hat z_i, \forall i$.
% The  following inequality will be satisfied,
% \begin{equation}\label{eq: reault1}
%     \mathcal L (\hat z^+, \hat\lambda^+) \leq \frac{1}{1+ \delta}\mathcal L (\hat z, \hat\lambda) + \frac{4}{1+ \delta} \mathcal O(N\Delta),
% \end{equation} 
% with Algorithm \ref{alg:ALADIN}, 
%     where %\[ \mathcal O(4 N \Delta) = \left(10\rho  M_z +  2 M_\lambda\right) N\Delta.\] 
%     \[\mathcal O(N \Delta) = 6\rho  M_z  N\Delta.\]
\end{theorem}
\textit{Proof.} See Appendix \ref{app1}. \hfill$\blacksquare$

% \begin{remark}
%     \todo{discuss what we see in Theorem 1 + what is the O(...), why it relies on Delta + what is the factor 4? (4 arises from Lemma 1 combining quantization errors for priman and dual?) + what happens if we increase / decrease Delta (i.e., if we increase we increase communication efficiency), if we decrease then we increase accuracy.} 
% \end{remark}

\begin{remark}\label{remark: tradeoff} 
    Focusing on \eqref{eq: reault1} of Theorem~\ref{the: convergence}, the term $\frac{4}{1+ \delta} \mathcal O(N\Delta)$ represents the quantization error introduced from Algorithm~\ref{alg:QuAS}. 
    As we will see later in Section~\ref{sec: numerical}, this error causes nodes to converge to a $\Delta$-dependent neighborhood of the optimal solution. 
    While decentralized approaches to progressively refine $\Delta$ can enhance solution precision \cite{rikos2023distributed}, this typically incurs higher communication overhead in terms of bits per message compromising our algorithm's communication efficiency.
    % Utilizing decentralized strategies for progressively refining $\Delta$ may enable nodes to approximate the optimal solution with even higher precision \cite{rikos2023distributed}. 
    % However, refining $\Delta$ incurs higher communication overhead in terms of bits per message compromising the communication efficiency of Algorithm~\ref{alg:QuAS}. 
    In contrast, employing quantizers with base-shifting capabilities \cite{rikos2024finite_bit_zooming} allows for maintaining high communication efficiency while still enabling nodes to approximate the optimal solution with greater precision.
    This latter strategy however, results in a trade-off as it may reduce the convergence speed of Algorithm~\ref{alg:QuAS}.
\end{remark}

\noindent
\textbf{Relaxing Strong Convexity.} %\todo{MODIFY THIS} 
While Theorem~\ref{the: convergence} relies on $\mu$-strong convexity (see Assumption~\ref{ass:2}), it is important to note that Algorithm~\ref{alg:ALADIN} can achieve a symmetric global linear convergence rate when the local cost function $f_i$ of each node $i \in \mathcal{V}$ is convex (and not strongly convex) provided it remains closed, proper, and $L$-smooth (i.e., satisfies Assumption~\ref{ass:3}). 
% While Theorem~\ref{the: convergence} relies on strong convexity (see Assumption~\ref{ass:2}), it is important to note that Algorithm~\ref{alg:ALADIN} can exhibit symmetric global linear convergence rate for the case where the local cost function $f_i$ of each node $i \in \mathcal{V}$ is satisfied with Assumption \ref{ass:2}.
% is convex (and not strongly convex) provided it remains closed, proper, smooth, and satisfies the Lipschitz condition (with constant $L_i$ as defined in \eqref{eq: lip}). 
We present this analysis in the following corollary. 

\begin{corollary}\label{coro1}
Let us consider a digraph $\mathcal G= \left( \mathcal V, \mathcal E\right)$. 
Each node $i \in \mathcal{V}$ has a local cost function $f_i$, %that is closed, proper, smooth, convex, and fulfills \eqref{eq: lip}. 
and Assumptions~\ref{ass:1} and \ref{ass:3} hold. 
%Also, Assumption~\ref{ass:1} holds. 
Each node $i \in \mathcal V$ executes Algorithm~\ref{alg:ALADIN} for solving the consensus optimization problem in~\eqref{eq: reformulate3} in a decentralized manner. 
During the operation of Algorithm~\ref{alg:ALADIN}, there always exists a $\delta>0$ such that
    \begin{equation}\label{eq: condition2}%\small
        \delta \mathcal{L}(\hat z^+, \hat \lambda^+)\leq 4 \sum_{i=1}^N \frac{1}{L_i}\left\|g_i-g_i^* \right\|^2.
    \end{equation} 
    As a result, we have that during the operation of Algorithm~\ref{alg:ALADIN} the inequality \eqref{eq: reault1} in Theorem~\ref{the: convergence} is satisfied for every node. 
    % From \eqref{eq: condition2}, equation \eqref{eq: reault1} is satisfied with  Algorithm \ref{alg:ALADIN}.
\end{corollary}
\textit{Proof.} See Appendix \ref{app: 2}. \hfill$\blacksquare$

\vspace{.1cm}

%Note that the error accumulation bound in Remark~\ref{remark error accumulation} remains valid also when considering the scenario described in Corollary~\ref{coro1}. 

% \noindent
% \textbf{Relaxing Strong Convexity and Lipschitz-continuity.}
% In the absence of Lipschitz-continuity and strong convexity, Algorithm~\ref{alg:ALADIN} can guarantee global convergence, but without providing any specific convergence rate \todo{(but can we prove we can do this?) -- maybe remove? because arxiv does not have a specific proof.} 

\section{Numerical Simulation}\label{sec: numerical}

% In this section, we present numerical simulations to demonstrate the operation of Algorithm~\ref{alg:ALADIN} and its potential advantages. 
% Additionally, we compare Algorithm~\ref{alg:ALADIN} against existing algorithms highlighting the improvements introduced. 

In this section, we present numerical simulations to demonstrate the operation of Algorithm~\ref{alg:ALADIN} and to highlight the improvements it offers over existing distributed optimization algorithms. 

We focus on a random digraph consisting of $20$ nodes. 
Each node $i$ is associated with a function $f_i(z) = \frac{1}{2} z^\top P_i z + p_i^\top z$ where $z\in \mathbb R^n$, $P_i\in \mathbb S_{++}^{n}$, and $p_i \in \mathbb R^n$ for each node $i \in \mathcal{V}$, with $n=20$. 
Furthermore, we have $\rho = 1$ and that Assumptions~\ref{ass:1} and~\ref{ass:2} hold.
For each node $i$, $P_i \succ 0$ was initialized as the square of a randomly generated symmetric matrix $A_i$, which guarantees that it is positive definite. 
Additionally, $q_i$ is set as the negative of the product of the transpose of $A_i$ and a randomly generated vector $b_i$ (i.e., it represents a linear term). 
For further details, please refer to \cite[Section~VI]{rikos2023asynchronous}.
We execute Algorithm~\ref{alg:ALADIN} and show how the nodes' decision variables convergence to the optimal solution for  $\Delta = 10^{-3},  10^{-4}$, and $10^{-5}$, respectively. 
We plot the error $\sum_{i=1}^N\left\|x_i^{[k]} - x_i^* \right\|$ where $x_i^*=z^*, \forall i \in \mathcal{V}$ represents the optimal solution of problem~\eqref{eq: reformulate3}. 

% In Fig.~\ref{fig: comparison}, we focus on a randomly generated digraph consisting of $20$ nodes. 
% Each node $i$ is given by a local cost function $f_i(z) = \frac{1}{2} z^\top P_i z + p_i^\top z$ where  $z\in \mathbb R^n$ and $P_i\in \mathbb S_{++}^{n}, p_i\in \mathbb R^n, \forall i$,  with $n=20$, such that the cost functions $f_i, \forall i$ are quadratic and convex. 
% This is consistent with the requirements of Assumption~\ref{ass:1} and~\ref{ass:2}. 
% Here the matrices $P_i$ and vectors $p_i$ are generated in the same way as \cite[Section~VI]{rikos2023asynchronous}. 
% We execute Algorithm \ref{alg:ALADIN} and show how the nodes' decision variables convergence to the optimal solution for $\Delta = 10^{-3},  10^{-4}, 10^{-5}$, respectively. 
% We plot the error $\sum_{i=1}^N\left\|x_i^{[k]} - x_i^* \right\|$ where $x_i^*=z^*, \forall i$ are the optimal solution of problem \eqref{eq: reformulate3}. 

\begin{figure}[h]
	\centering
\includegraphics[width=0.5\textwidth,height=0.29\textheight]{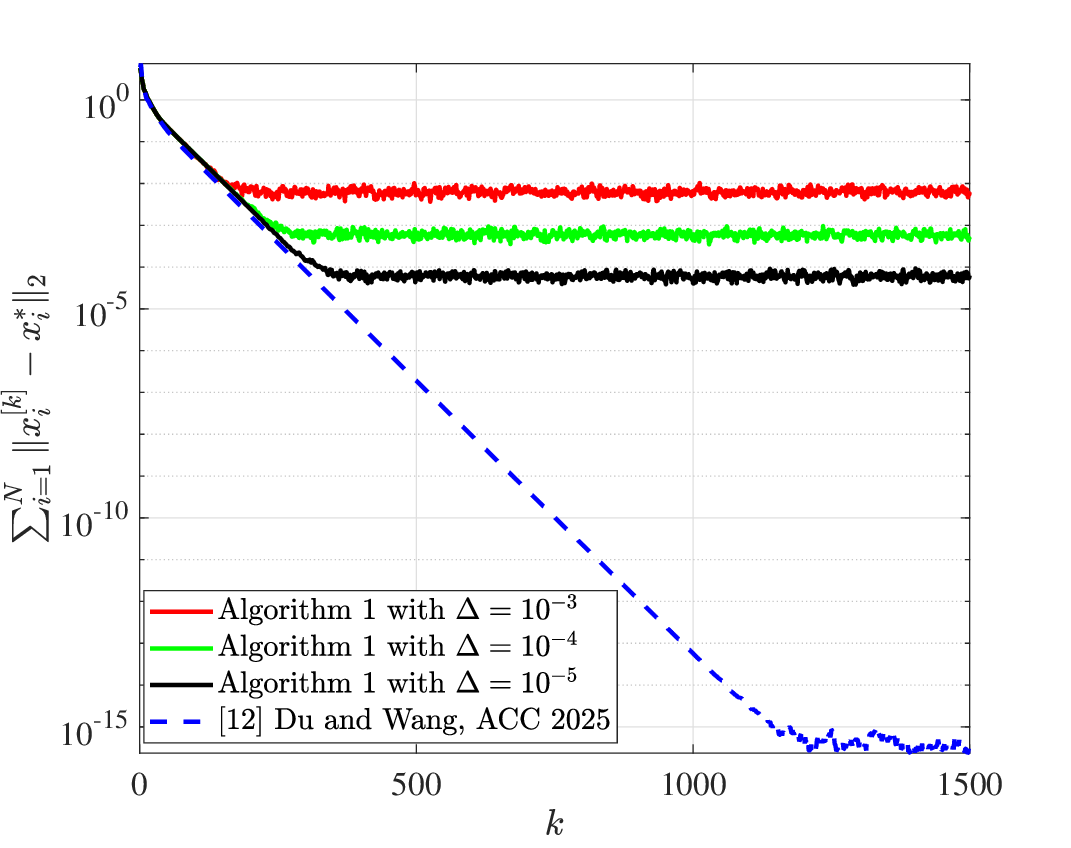}
	\caption{Comparison of Algorithm \ref{alg:ALADIN} 
    with RC-ALADIN \cite{Du2023} over a directed graph with quantization level $\Delta = 10^{-3},  10^{-4}$ and $10^{-5}$.}
	\label{fig: comparison}
\end{figure}

In Fig.~~\ref{fig: comparison} we can see that Algorithm~\ref{alg:ALADIN} converges to the optimal solution, achieving an approximation precision that is directly influenced by the quantization level $\Delta$. 
Specifically, a smaller $\Delta$ enables nodes to approximate the optimal solution with greater accuracy, a behavior consistent with the theoretical results presented in Theorem~\ref{the: convergence}. 
The oscillatory behavior observed in the convergence is attributed to the nonlinearities introduced by quantized communication, which affect the values of parameters such as $x_i$, $g_i$, and $\hat \lambda_i$. 
While Algorithm~\ref{alg:ALADIN} demonstrates comparable performance to the approach in \cite{Du2023} up to the neighborhood of the optimal solution, \cite{Du2023} relies on a distributed framework with a server node and real-valued message exchanges among nodes.
Thus, the approach in \cite{Du2023} introduces scalability challenges, suffers from single-point-of-failure risks, limits its applicability to bandwidth-constrained environments and compromises resource efficiency. 
In contrast, Algorithm~\ref{alg:ALADIN} offers a significant advantages. 
It provides comparable performance with algorithms in the literature that enable nodes to exchange real-valued messages (see \cite{Du2023}), while exhibiting efficient (quantized) communication among nodes {\color{black} {(a larger $\Delta$ implies a reduced communication bandwidth requirement among nodes)}} thereby reducing data transmission. 
Moreover, Algorithm~\ref{alg:ALADIN} operates in a fully decentralized manner, eliminating the need for a server node. {\color{black}{A numerical example in \cite[Section 5.1]{Houska2021} investigates the role of the penalty parameter $\rho$ in Algorithm~\ref{alg:ALADIN}, applied to a consistent convex problem under Assumptions~\ref{ass:2} and \ref{ass:3}. The example demonstrates that $\rho$ influences the observed linear convergence rate. Nevertheless, a comprehensive theoretical characterization of this dependence for the full class of problems satisfying the same assumptions remains an open question.}}

\section{Conclusions and Future Directions}\label{sec: conclusion}

In this paper, we presented a novel decentralized optimization algorithm named QuDRC-ALADIN. 
Our algorithm enables computation of the optimal solution in a fully decentralized manner over directed communication networks, while ensuring efficient quantized communication among nodes. 
We analyzed our algorithm's operation and established its linear convergence to a neighborhood of the optimal solution that depends on the utilized quantization level. 
Finally, we presented numerical simulations validating our algorithm’s performance, and highlighted its advantages compared to existing algorithms in the literature. 

Future work will focus on extending the algorithm to non-convex optimization problems, and analyzing asynchronous communication protocols for enhancing robustness and scalability in dynamic environments. 

\appendices
\section{Proof of Theorem \ref{the: convergence}}\label{app1}

The analysis is inspired by \cite{Du2023} and \cite{jiang2021asynchronous}. Before starting the proof, let us recall the following equality 
\begin{equation}\label{eq: well-known equality}\small
2\left(\alpha-\gamma \right)^\top \left(\alpha-\beta \right) =
    \left\|\alpha-\gamma \right\|^2 - \left\|\beta-\gamma \right\|^2+\left\|\alpha-\beta \right\|^2,
\end{equation}
% \begin{equation}\label{eq: well-known equality}%\small
% \begin{split}
%     &2\left(\alpha-\gamma \right)^\top \left(\alpha-\beta \right) \\
%     =& \left\|\alpha-\gamma \right\|^2 - \left\|\beta-\gamma \right\|^2+\left\|\alpha-\beta \right\|^2,\; \forall \alpha, \beta, \gamma \in \mathbb R^n.
% \end{split}
% \end{equation}
for $\alpha, \beta, \gamma \in \mathbb R^n$. 
Focusing on equation \eqref{eq: mu}, the following symmetric inequalities can be obtained 
\begin{equation}\small
\label{eq: strongly convex}\left\{
\begin{split}
    & f_i(x_\alpha)+\nabla f_i(x_\alpha)^\top (x_\beta-x_\alpha) +\frac{\mu_i}{2}\|x_\beta-x_\alpha\|^2\leq f_i(x_\beta)\\
    & f_i(x_\beta)+\nabla f_i(x_\beta)^\top (x_\alpha-x_\beta) +\frac{\mu_i}{2}\|x_\alpha-x_\beta\|^2\leq f_i(x_\alpha).
\end{split}\right.
\end{equation}
By adding up the two inequalities of \eqref{eq: strongly convex}, the following equation can be derived 
\begin{equation}\label{after_strong_convex}\small
    \mu_i \left\| x_\alpha -x_\beta\right\|^2\leq\left(\nabla f_i(x_\alpha) -\nabla f_i(x_\beta)\right)^\top \left(x_\alpha -x_\beta \right). 
\end{equation}
By setting $x_\alpha= x_i^+$ and $x_\beta = z^*$ in \eqref{after_strong_convex}, we have
\begin{equation}\label{eq: strongly convex2}\small
    \mu_i \left\| x_i^+ -z^*\right\|^2\leq\left(\nabla f_i(x_i^+) -\nabla f_i(z^*)\right)^\top \left(x_i^+ -z^* \right).
\end{equation}
By summing the $N$ inequalities in \eqref{eq: strongly convex2}, we have 
\begin{equation}\label{eq: strongly convex3}\small
   \sum_{i=1}^N \mu_i \left\| x_i^+ -z^*\right\|^2 \leq \sum_{i=1}^N\left (\nabla f_i(x_i^+) -\nabla f_i(z^*)\right)^\top \left(x_i^+ -z^* \right).
\end{equation}
Note that due to the simple smoothness of $f_i, \forall i$ in Assumption \ref{ass:2},  we have $g_i=\nabla f_i$. 
Thus, the right hand side of \eqref{eq: strongly convex3} can be represented as,
\begin{align}\label{eq: simplify}\small
    \begin{split}
        % &\sum_{i=1}^N\left (\nabla f_i(x_i^+) -\nabla f_i(z^*)\right)^\top \left(x_i^+ -z^* \right)\\
        % {=}
        & \sum_{i=1}^N\left (g_i -g_i^*\right)^\top \left(x_i^+ -z^* \right)\\
        \overset{\eqref{eq: new gradient}}{=}& \sum_{i=1}^N\left (\rho\left(\hat z_i-x_i^+\right)-\hat\lambda_i +\lambda_i^*\right)^\top \left(x_i^+ -z^* \right)\\
        \overset{\eqref{eq: lemma}}{=}& \sum_{i=1}^N\left (\rho\left(
        \frac{\hat z_i-\hat z^+_i}{2}
        -\frac{\hat \lambda_i^+-\hat\lambda_i}{2\rho}\right)-\hat \lambda_i +\lambda_i^*\right)^\top \\
        &\hspace{2.5mm}\quad\left(\frac{\hat \lambda_i^+-\hat\lambda_i}{2\rho}+\frac{\hat z^+_i+\hat z_i}{2} -z^* \right).
    \end{split}
\end{align}
We now split the right-hand side of \eqref{eq: simplify} into the following components (a), (b) (c), (d), (e). 
Then, we analyze each component separately. 
 \vspace{-2mm}
     %\vspace{-1mm}\begin{subequations}
    {\small \begin{align}
   (a) &\quad         \frac{1}{4}\sum_{i=1}^N \left(\hat z_i - \hat z_i^+ \right)^\top \left( \hat\lambda_i^+ - \hat\lambda_i\right)\notag\\
          & = \frac{1}{4}\sum_{i=1}^N \left(\left( z_i-z^*\right) - \left( z_i^+-z^* \right)\right.\label{eq: 1}\\
           &\hspace{2mm}\qquad\left.+\left( \hat z_i-z_i \right) -\left( \hat z_i^+-z_i^+ \right) \right)^\top \notag\\
           &\;\qquad\left(\left( \lambda_i-\lambda_i^*\right) - \left( \lambda_i^+-\lambda_i^* \right)\right.\notag\\
           &\hspace{2mm}\qquad\left.+\left( \hat \lambda_i-\lambda_i \right) -\left( \hat \lambda_i^+-\lambda_i^+ \right) \right)\notag\\
          &  \overset{\eqref{eq: error},\eqref{eq: lam}}{\leq} 2N\rho \left(M_z\Delta+2\Delta^2 \right). \notag\\
         (b) &\quad  \sum_{i=1}^N \left( \frac{\hat \lambda_i^+ - \hat \lambda_i}{2\rho}\right)^\top\left( \frac{\hat \lambda_i -\hat \lambda_i^+}{2} +\lambda_i^* -\hat\lambda_i\right)\label{eq: 23}\\
       &     \overset{\eqref{eq: well-known equality}}{=} \frac{1}{4\rho}\sum_{i=1}^N \left(\|\hat \lambda_i-\lambda_i^*\|^2 -\|\hat \lambda_i^+-\lambda_i^*\|^2\right).\notag\\
            % \overset{\eqref{eq: error},\eqref{eq: lam}}{\leq}&\frac{1}{4\rho}\sum_{i=1}^N \hspace{-1mm}\left(\| \lambda_i-\lambda_i^*\|^2 \hspace{-1mm}-\hspace{-1mm}\| \lambda_i^+-\lambda_i^*\|^2\right)\hspace{-1mm} \\
            % &%+\frac{N}{2}\hspace{-1mm}\left(2M_\lambda \Delta +  \Delta^2 \right).
            % +2N\left(M_\lambda \Delta +  \rho \Delta^2 \right).
          (c)& \quad \frac{\rho}{2}\sum_{i=1}^N \left( \frac{\hat z_i^++\hat z_i}{2}-z^*\right)^\top\left( \hat z_i - \hat z_i^+ \right)\label{eq: 4+7}\\
         &   \overset{\eqref{eq: well-known equality}}{=} \frac{\rho}{4} \sum_{i=1}^N \left( \left\|\hat z_i - z^* \right\|^2 -\left\|\hat z_i^+ - z^* \right\|^2 \right).\notag\\
            % \overset{\eqref{eq: error}}{\leq} &
            % \frac{\rho }{4}\sum_{i=1}^N \left( \left\| z_i - z^* \right\|^2 -\left\|z_i^+ - z^* \right\|^2 \right)\\
            % &%+ \frac{\rho N}{2}\left(2M_z \Delta + \Delta^2 \right).
            % +2N\rho \left(M_z \Delta +  \Delta^2 \right).
   (d)&   \quad      \sum_{i=1}^N \left(\frac{\hat z_i^++\hat z_i -2z^*}{2} \right)^\top \left(\frac{\hat \lambda_i - \hat\lambda_i^+}{2} \right)\label{eq: 5+8}\\
         &   \overset{\eqref{eq: error},\eqref{eq: lam}}{\leq}  2N\rho \left(M_z\Delta+2\Delta^2 \right).\notag\\
         % \end{align}}
         %  {\small \begin{align}
       (e)&\quad     \sum_{i=1}^N \left(\frac{\hat z_i^++\hat z_i -2z^*}{2} \right)^\top \left(\lambda_i^* -\hat \lambda_i \right)\label{eq: 6+9}\\
           &\overset{\eqref{eq: error},\eqref{eq: lam}}{\leq}  2N\rho \left(M_z\Delta+2\Delta^2 \right).\notag
    \end{align}}
  %  \end{subequations}
    % \begin{equation}\label{eq: 5+8}\small
    %     \begin{split}
    %         &\sum_{i=1}^N \left(\frac{\hat z_i^++\hat z_i -2z^*}{2} \right)^\top \left(\frac{\hat \lambda_i - \hat\lambda_i^+}{2} \right)\\
    %         =&\sum_{i=1}^N \left(\frac{ ( \hat z_i^+-z^*)   +(\hat z_i -z^*)}{2} \right)^\top \\
    %         &\hspace{2mm}\quad\left(\frac{\left(\hat \lambda_i-\lambda_i^* \right) -  \left(\hat\lambda_i^+-\lambda_i^*\right)}{2} \right)\\
    %         \overset{\eqref{eq: error},\eqref{eq: lam}}{\leq} & 2N\rho \left(M_z\Delta+2\Delta^2 \right).
    %     \end{split}
    % \end{equation}
     %\vspace{-1mm}

    % \begin{equation}\label{eq: 6+9}\small
    %     \begin{split}
    %         &\sum_{i=1}^N \left(\frac{\hat z_i^++\hat z_i -2z^*}{2} \right)^\top \left(\lambda_i^* -\hat \lambda_i \right)\\
    %         =&\sum_{i=1}^N \left(\frac{(\hat z_i^+-z^*)   +(\hat z_i -z^*)}{2} \right)^\top \left(\lambda_i^* -\hat \lambda_i \right)\\
    %        \overset{\eqref{eq: error},\eqref{eq: lam}}{\leq}&  2N\rho \left(M_z\Delta+2\Delta^2 \right).
    %     \end{split}
    % \end{equation}

Note here that the analysis in \eqref{eq: 5+8} and \eqref{eq: 6+9} followed the same methodology as the one in \eqref{eq: 1}, but we omitted the details for space considerations. 
% \todo{CHECK THIS AGAIN -- APOSTOLOS}
% \todo{Following the same derivation as \eqref{eq: 1}, equation \eqref{eq: 5+8} and \eqref{eq: 6+9} are upper-bounded by $2N\rho \left(M_z\Delta+2\Delta^2 \right)$,  and the details are therefore omitted here.} 
Considering now that $\hat z= \hat z_i$ for every node $i \in \mathcal{V}$, and summing \eqref{eq: 1} -- \eqref{eq: 6+9}, we have that \eqref{eq: strongly convex3} becomes 
% With equation \eqref{eq: strongly convex3} -- \eqref{eq: 6+9}, and $\hat z= \hat z_i, \forall i$, we can show 
\vspace{-1mm}
{\small\begin{align}\label{eq: end}
    \sum_{i=1}^N \mu_i \left\| x_i^+ -z^*\right\|^2 \leq \frac{1}{4}\mathcal L \left(\hat z, \hat\lambda\right)-\frac{1}{4}\mathcal L \left(\hat z^+, \hat\lambda^+\right)+ \mathcal O(N \Delta)
\end{align}}
% \begin{equation}\label{eq: end}\small
%     \begin{split}
%         &\sum_{i=1}^N \mu_i \left\| x_i^+ -z^*\right\|^2 \\\leq& \frac{1}{4}\mathcal L \left(\hat z, \hat\lambda\right)-\frac{1}{4}\mathcal L \left(\hat z^+, \hat\lambda^+\right)+ \mathcal O(N \Delta),
%     \end{split}
% \end{equation}
 %\[\mathcal O(N \Delta) = \left(M_\lambda+ 4\rho  M_z \right) N\Delta.\]
where $\mathcal O(N \Delta) = 6\rho  M_z  N\Delta$. 
Note that in \eqref{eq: end} we omitted the higher-order terms of $\Delta$, as they are significantly smaller than the first-order terms and can be considered negligible. 

%Finally, combining \eqref{eq: end} with \eqref{eq: condition}, we have that \eqref{eq: reault1} holds. 

% Finally, combining \eqref{eq: end} with \eqref{eq: condition}, we can obtain the following equation, $\frac{\delta}{4}\mathcal L \left(\hat z^+, \hat\lambda^+ \right)\leq \frac{1}{4}\mathcal L\left(\hat z, \hat \lambda \right) - \frac{1}{4}\mathcal L\left(\hat z^+, \hat \lambda^+ \right)+ \mathcal O(N \Delta),$
\begin{equation}\label{eq: range}\small
    \frac{\delta}{4}\mathcal L \left(\hat z^+, \hat\lambda^+ \right)\leq \frac{1}{4}\mathcal L\left(\hat z, \hat \lambda \right) - \frac{1}{4}\mathcal L\left(\hat z^+, \hat \lambda^+ \right)+ \mathcal O(N \Delta),
\end{equation}
which is equivalent to \eqref{eq: reault1}.
This concludes the proof of our theorem.

\section{Proof of Corollary \ref{coro1}}\label{app: 2}

Focusing on \eqref{eq: lip}, considering that $g_i=\nabla f_i$, setting $x_\alpha= x_i^+$, $x_\beta = z^*$, and summing the $N$ inequalities we have 
\vspace{-1mm}
\begin{equation}\label{eq: last}\small
   \sum_{i=1}^N \frac{1}{L_i}\left\| g_i -g_i^*\right\|^2\leq\sum_{i=1}^N\left ( g_i -g_i^*\right)^\top \left(x_i^+ -z^* \right).
\end{equation}

Note that, the right-hand side of \eqref{eq: last} is the same as the right-hand side of \eqref{eq: strongly convex3}. 
From our analysis in Appendix~\ref{app1}, this means that the right-hand side of \eqref{eq: last} is upper bounded by the right-hand side of \eqref{eq: end}, therefore the following equation is satisfied 
\vspace{-1mm}
\begin{equation}\label{eq: end2}\small
    \begin{split}
        \sum_{i=1}^N \frac{1}{L_i}\left\| g_i -g_i^*\right\|^2 \leq \frac{1}{4}\mathcal L \left(\hat z, \hat\lambda\right)-\frac{1}{4}\mathcal L \left(\hat z^+, \hat\lambda^+\right)+ \mathcal O(N \Delta).
    \end{split}
\end{equation} 
Finally, combining inequality \eqref{eq: condition2} and inequality \eqref{eq: end2}, we can then prove inequality \eqref{eq: reault1}.
\bibliographystyle{IEEEtran}
		\bibliography{paper}	

\begin{thebibliography}{10}
\providecommand{\url}[1]{#1}
\csname url@samestyle\endcsname
\providecommand{\newblock}{\relax}
\providecommand{\bibinfo}[2]{#2}
\providecommand{\BIBentrySTDinterwordspacing}{\spaceskip=0pt\relax}
\providecommand{\BIBentryALTinterwordstretchfactor}{4}
\providecommand{\BIBentryALTinterwordspacing}{\spaceskip=\fontdimen2\font plus
\BIBentryALTinterwordstretchfactor\fontdimen3\font minus \fontdimen4\font\relax}
\providecommand{\BIBforeignlanguage}[2]{{%
\expandafter\ifx\csname l@#1\endcsname\relax
\typeout{** WARNING: IEEEtran.bst: No hyphenation pattern has been}%
\typeout{** loaded for the language `#1'. Using the pattern for}%
\typeout{** the default language instead.}%
\else
\language=\csname l@#1\endcsname
\fi
#2}}
\providecommand{\BIBdecl}{\relax}
\BIBdecl

\bibitem{gosta2024decentralized}
G.~Stomberg, A.~Engelmann, M.~Diehl, and T.~Faulwasser, ``Decentralized real-time iterations for distributed nonlinear model predictive control,'' \emph{arXiv preprint arXiv:2401.14898}, 2024.

\bibitem{decentralizedlearningsurvey}
L.~Yuan, Z.~Wang, L.~Sun, P.~S. Yu, and C.~G. Brinton, ``Decentralized federated learning: A survey and perspective,'' \emph{IEEE Internet of Things Journal}, vol.~11, no.~21, pp. 34\,617--34\,638, 2024.

\bibitem{liang2024decentralized}
D.~Liang, S.~Su, L.~Zeng, and H.-D. Chiang, ``Decentralized method for nonconvex robust static state estimation of integrated electricity-gas systems,'' \emph{CSEE Journal of Power and Energy Systems}, pp. 1--13, 2024.

\bibitem{nedic2009distributed}
A.~Nedic and A.~Ozdaglar, ``Distributed subgradient methods for multi-agent optimization,'' \emph{IEEE Transactions on Automatic Control}, vol.~54, no.~1, pp. 48--61, 2009.

\bibitem{wang2018giant}
S.~Wang, F.~Roosta, P.~Xu, and M.~W. Mahoney, ``{GIANT}: Globally improved approximate newton method for distributed optimization,'' in \emph{Advances in Neural Information Processing Systems}, vol.~31, 2018.

\bibitem{boyd2011distributed}
S.~Boyd, N.~Parikh, E.~Chu, B.~Peleato, J.~Eckstein \emph{et~al.}, ``Distributed optimization and statistical learning via the alternating direction method of multipliers,'' \emph{Foundations and Trends{\textregistered} in Machine learning}, vol.~3, no.~1, pp. 1--122, 2011.

\bibitem{Houska2016}
B.~Houska, J.~Frasch, and M.~Diehl, ``An augmented {L}agrangian based algorithm for distributed {nonconvex} optimization,'' \emph{SIAM Journal on Optimization}, vol.~26, no.~2, pp. 1101--1127, 2016.

\bibitem{ling2015dlm}
Q.~Ling, W.~Shi, G.~Wu, and A.~Ribeiro, ``{DLM}: Decentralized linearized alternating direction method of multipliers,'' \emph{IEEE Transactions on Signal Processing}, vol.~63, no.~15, pp. 4051--4064, 2015.

\bibitem{dantzig1960decomposition}
G.~B. Dantzig and P.~Wolfe, ``Decomposition principle for linear programs,'' \emph{Operations research}, vol.~8, no.~1, pp. 101--111, 1960.

\bibitem{he20121}
B.~He and X.~Yuan, ``On the {O}(1/n) convergence rate of the douglas--rachford alternating direction method,'' \emph{SIAM Journal on Numerical Analysis}, vol.~50, no.~2, pp. 700--709, 2012.

\bibitem{Hong2016}
M.~Hong, Z.-Q. Luo, and M.~Razaviyayn, ``Convergence analysis of alternating direction method of multipliers for a family of nonconvex problems,'' \emph{SIAM Journal on Optimization}, vol.~26, no.~1, pp. 337--364, 2016.

\bibitem{Du2023}
X.~Du and J.~Wang, ``Distributed consensus optimization with consensus {ALADIN},'' in \emph{American Control Conference}, 2025 (accepted for publication).

\bibitem{Houska2021}
B.~Houska and Y.~Jiang, ``Distributed optimization and control with aladin,'' \emph{Recent Advances in Model Predictive Control: Theory, Algorithms, and Applications}, pp. 135--163, 2021.

\bibitem{Engelmann2019}
A.~Engelmann, Y.~Jiang, T.~M{\"u}hlpfordt, B.~Houska, and T.~Faulwasser, ``Toward distributed {OPF} using {ALADIN},'' \emph{IEEE Transactions on Power Systems}, vol.~34, no.~1, pp. 584--594, 2019.

\bibitem{Du2019}
X.~Du, A.~Engelmann, Y.~Jiang, T.~Faulwasser, and B.~Houska, ``Distributed state estimation for {AC} power systems using {G}auss-{N}ewton {ALADIN},'' in \emph{IEEE Conference on Decision and Control}, 2019, pp. 1919--1924.

\bibitem{falsone2023augmented}
A.~Falsone and M.~Prandini, ``Augmented {Lagrangian} tracking for distributed optimization with equality and inequality coupling constraints,'' \emph{Automatica}, vol. 157, p. 111269, 2023.

\bibitem{jiang2021asynchronous}
W.~Jiang, A.~Grammenos, E.~Kalyvianaki, and T.~Charalambous, ``An asynchronous approximate distributed alternating direction method of multipliers in digraphs,'' in \emph{IEEE Conference on Decision and Control}, 2021, pp. 3406--3413.

\bibitem{engelmann2020decomposition}
A.~Engelmann, Y.~Jiang, B.~Houska, and T.~Faulwasser, ``Decomposition of nonconvex optimization via bi-level distributed {ALADIN},'' \emph{IEEE Transactions on Control of Network Systems}, vol.~7, no.~4, pp. 1848--1858, 2020.

\bibitem{Du2023B}
X.~Du, J.~Wang, X.~Zhou, and Y.~Mao, ``A bi-level globalization strategy for non-convex consensus {ADMM} and {ALADIN},'' \emph{arXiv preprint arXiv:2309.02660}, 2023.

\bibitem{zhu2016quantized}
S.~Zhu, M.~Hong, and B.~Chen, ``Quantized consensus {ADMM} for multi-agent distributed optimization,'' in \emph{IEEE International Conference on Acoustics, Speech and Signal Processing}, 2016, pp. 4134--4138.

\bibitem{elgabli2020q}
A.~Elgabli, J.~Park, A.~S. Bedi, C.~B. Issaid, M.~Bennis, and V.~Aggarwal, ``{Q-GADMM}: Quantized group {ADMM} for communication efficient decentralized machine learning,'' \emph{IEEE Transactions on Communications}, vol.~69, no.~1, pp. 164--181, 2020.

\bibitem{rikos2023asynchronous}
A.~I. Rikos, W.~Jiang, T.~Charalambous, and K.~H. Johansson, ``Asynchronous distributed optimization via {ADMM} with efficient communication,'' in \emph{IEEE Conference on Decision and Control}, 2023, pp. 7002--7008.

\bibitem{rikos2025distributed}
------, ``Distributed optimization with efficient communication, event-triggered solution enhancement, and operation stopping,'' \emph{arXiv preprint arXiv:2504.16477}, 2025.

\bibitem{proakis2002communication}
J.~G. Proakis and M.~Salehi, \emph{Communication Systems Engineering}, 2nd~ed.\hskip 1em plus 0.5em minus 0.4em\relax Upper Saddle River, N.J.: Prentice Hall, 2002.

\bibitem{2019:Wei_Johansson}
J.~Wei, X.~Yi, H.~Sandberg, and K.~H. Johansson, ``Nonlinear consensus protocols with applications to quantized communication and actuation,'' \emph{IEEE Transactions on Control of Network Systems}, vol.~6, no.~2, pp. 598--608, 2019.

\bibitem{Du2023_Arxiv}
X.~Du and J.~Wang, ``Consensus {ALADIN}: A framework for distributed optimization and its application in federated learning,'' \emph{arXiv preprint arXiv:2306.05662}, 2023.

\bibitem{OLIVA201620}
G.~Oliva, R.~Setola, and C.~N. Hadjicostis, ``Distributed finite-time calculation of node eccentricities, graph radius and graph diameter,'' \emph{Systems $\&$ Control Letters}, vol.~92, pp. 20--27, 2016.

\bibitem{deplano2025optimization}
D.~Deplano, N.~Bastianello, M.~Franceschelli, and K.~H. Johansson, ``Optimization and learning in open multi-agent systems,'' \emph{arXiv preprint arXiv:2501.16847}, 2025.

\bibitem{beck2017first}
A.~Beck, \emph{First-order methods in optimization}.\hskip 1em plus 0.5em minus 0.4em\relax SIAM, 2017.

\bibitem{Rawlings2017}
J.~B. Rawlings, D.~Q. Mayne, and M.~Diehl, \emph{Model Predictive Control: Theory, Computation, and Design, 2nd Edition}.\hskip 1em plus 0.5em minus 0.4em\relax Nob Hill Publishing, 2017.

\bibitem{boyd2004convex}
S.~P. Boyd and L.~Vandenberghe, \emph{Convex optimization}.\hskip 1em plus 0.5em minus 0.4em\relax Cambridge university press, 2004.

\bibitem{rikos2022non}
A.~I. Rikos, C.~N. Hadjicostis, and K.~H. Johansson, ``Non-oscillating quantized average consensus over dynamic directed topologies,'' \emph{Automatica}, vol. 146, p. 110621, 2022.

\bibitem{rikos2024distributed}
A.~I. Rikos, A.~Grammenos, E.~Kalyvianaki, C.~N. Hadjicostis, T.~Charalambous, and K.~H. Johansson, ``Distributed optimization for quadratic cost functions with quantized communication and finite-time convergence,'' \emph{IEEE Transactions on Control of Network Systems}, vol.~12, no.~1, pp. 930--942, 2025.

\bibitem{rikos2023distributed}
A.~I. Rikos, W.~Jiang, T.~Charalambous, and K.~H. Johansson, ``Distributed optimization via gradient descent with event-triggered zooming over quantized communication,'' in \emph{IEEE Conference on Decision and Control}, 2023, pp. 6321--6327.

\bibitem{rikos2023distributed2}
------, ``Distributed optimization with gradient descent and quantized communication,'' \emph{IFAC-PapersOnLine}, vol.~56, no.~2, pp. 5900--5906, 2023.

\bibitem{rikos2024finite_bit_zooming}
------, ``Distributed optimization with finite bit adaptive quantization for efficient communication and precision enhancement,'' in \emph{IEEE Conference on Decision and Control}, 2024, pp. 2531--2537.

\end{thebibliography}

\end{document}